\documentclass[10pt]{amsart}
\usepackage{latexsym}
\usepackage{amsmath}
\usepackage{amssymb}
\usepackage{mathrsfs}
\usepackage{calrsfs}
\usepackage{bbold}

\def\card{\mathop{\fam0 card}}
 \def\Orth{\mathop{\fam0 Orth}}
 \def\ZFC{\mathop{\fam0 ZFC}}

 \def\qca{\mathop{\fam0 qca}}
 \def\Hom{\mathop{\fam0 Hom}}
 \def\upwardarrow{\mathord{
  \hbox to 5pt{\hss$\vcenter{\hbox to 2.4pt{\hss$\mathchar"222$\hss}\hrule}\hss$}
}}
\def\downwardarrow{\mathord{
  \hbox to 5pt{\hss$\vcenter{\hrule\hbox to 2.4pt{\hss$\mathchar"223$\hss}}\hss$}
}}

 \def\mix{\mathop{\fam0 mix}\nolimits}

 \def\[{\mathopen{\kern1pt\/ \vrule height7.6pt depth2.2pt width1pt\kern1.5pt}}
\def\]{\mathclose{\kern1.5pt\vrule height7.6pt depth2.2pt width1pt\kern1pt}}
\def\Bnorml{\mathopen{\kern1pt\/ \vrule height17pt depth10.2pt width1pt\kern1.5pt}}
\def\Bnormr{\mathclose{\kern1.5pt \vrule height17pt depth10.2pt width1pt\kern1pt}}
\def\bnorml{\mathopen{\kern1pt\/ \vrule height13pt depth7.2pt width1pt\kern1.5pt}}
\def\bnormr{\mathclose{\kern1.5pt \vrule height13pt depth7.2pt width1pt\kern1pt}}

\newcommand{\Theorem}[1]{\smallskip\textbf{{Theorem#1}.\/ }\sl}
\newcommand{\theorem}[1]{{\textbf{Theorem#1.\/~}}\sl}
\newcommand{\Corollary}[1]{\smallskip\textbf{{Corollary#1}.\/~}\sl}

\newcommand{\Proclaim}[1]{\smallskip{\bf#1}~\sl}
\newcommand{\proclaim}[1]{{\bf#1}~\sl}
\newcommand{\Endproc}{\rm}

 \def\subsec#1{\smallskip\textbf{#1}}
\def\osum{\mathop{o\text{\/-}\!\sum}}  

\def\olim{\mathop{o\text{-}{\fam0 lim}}} 

\def\beginproof{\par\mbox{$\vartriangleleft$}}
\def\endproof{\text{$\vartriangleright$}\smallskip}

\begin{document}

 \setcounter{footnote}{0}
 \setcounter{equation}{0}


 \title
[Some Aspects of Boolean Valued Analysis]
{Some Aspects of Boolean Valued Analysis}

 \author{A.~G.~Kusraev and S.~S.~Kutateladze}

 \begin{abstract}
This is a~survey of some recent applications of Boolean valued analysis
 to operator theory and harmonic analysis.
 Under consideration are pseudoembedding operators, the noncommutative Wickstead problem,
  the Radon--Nikod\'ym Theorem for $JB$-algebras, and the Bochner Theorem for lattice-valued
  positive definite mappings on locally compact groups.
\end{abstract}

\date{February 14, 2015}

\maketitle

 \section*{1.~Introduction}

 We survey here some aspects of Boolean valued analysis that concern operator theory.
 The term {\it Boolean valued analysis\/} signifies the technique
 of studying properties of an arbitrary mathematical
 object by  comparison between its representations
 in two different set-theoretic models whose construction
 utilizes principally distinct Boolean algebras.
 As these models, the classical Cantorian paradise in
 the shape of the von Neumann universe $\mathbb{V}$ and
 a specially-trimmed Boolean valued universe
 $\mathbb{V}^{(\mathbb{B})}$ are usually taken. Comparison
 analysis is carried out by  some interplay between
 $\mathbb{V}$ and $\mathbb{V}^{(\mathbb{B})}$.

 Boolean valued analysis not only is tied up with many topological and geometrical ideas but also provides a technology for expanding the content of the already available theorems. Each theorem, proven by the classical means, possesses some 
new unobvious content that relates
 to ``variable sets.'' A general scheme of the method is as follows; see \cite{IBA, BA_ST}. Assume that
 $\mathbf{X}\subset\mathbb{V}$ and $\mathbb{X}\subset\mathbb{V}^{(\mathbb{B})}$ are two classes of mathematical objects. Suppose that we are able to prove
 \smallskip

 The {\it Boolean Valued Representation:}~Every $X\in\mathbf{X}$ embeds into a Boolean valued model, becoming an object $\mathcal{X}\in\mathbb{X}$ within $\mathbb{V}^{(\mathbb{B})}$.
 \smallskip

 The {\it Boolean Valued Transfer Principle}  tells us then that every theorem about $\mathcal{X}$ within Zermelo--Fraenkel set theory has its counter\-part for the original object $X$ interpreted as a Boolean valued object $\mathcal{X}$.
 \smallskip

  The {\it Boolean Valued Machinery} enables us to perform some translation of theorems from $\mathcal{X}\in\mathbb{V}^{(\mathbb{B})}$ to $X\in\mathbb{V}$ by using the appropriate general operations (ascending--descending) and the principles of Boolean valued analysis.

  Everywhere below $\mathcal{R}$ stands for the reals within $\mathbb{V}^{(\mathbb{B})}$. The Gordon Theorem  states that the descent $\mathcal{R}{\downarrow}$ which is an algebraic system in $\mathbb{V}$ is a universally complete vector lattice; see \cite{IBA, BA_ST}.  Moreover, there exists an~isomorphism~$\chi$ of~$\mathbb{B}$ onto the Boolean algebra $\mathbb{P}(\mathrsfs R{\downarrow})$ such that
  $$
 \gathered
 \chi (b) x=\chi (b) y\Longleftrightarrow b\le [\![\,x=y\,]\!],
 \\
 \chi (b) x\le \chi (b) y\Longleftrightarrow b\le [\![\,x\le
 y\,]\!]
 \endgathered
 \eqno(\mathbb{G})
 $$
  for all~$x, y\in \mathrsfs R{\downarrow}$ and $b\in\mathbb{B}$. The {\it restricted descent} $\mathcal{R}{\Downarrow}$ of $\mathcal{R}$ is the part of  $\mathcal{R}{\downarrow}$ consisting of elements $x\in\mathcal{R}{\downarrow}$ with $|x|\leq C\mathbf{1}$ for some $C\in\mathbb{R}$, where $\mathbf{1}$ is an order unit in $\mathcal{R}{\downarrow}$.

 By a~{\it vector lattice}  throughout the~sequel we mean an Archimedean real vector lattice. We denote the~Boolean
 algebras of all bands and all band projections in a~vector lattice~$X$ respectively by $\mathbb{B}(X)$ and $\mathbb{P}(X)$   and we let $\mathcal{Z}(X)$ and $\Orth(X)$ stand for the {\it ideal center} of $X$ and the $f$-{\it algebra of orthomorphisms} on $X$, respectively. The~universal completion $X^\mathbb{u}$ of a vector lattice $X$ is always considered as a semiprime $f$-algebra whose multiplication is uniquely determined by fixing an order unit as a ring unit. The space of all order bounded linear operators from $X$ to $Y$ is denoted by $L^\sim(X,Y)$. The Riesz--Kan\-to\-rovich Theorem tells us that if $Y$ is a Dedekind complete vector lattice then so is $L^\sim(X,Y)$.

 A linear operator $T$ from $X$ to $Y$ is a \textit{lattice homomorphism\/} whenever $T$ preserves
 lattice operations; i.e., $T(x_1\vee x_2)=Tx_1\vee Tx_2$ (and so
 $T(x_1\wedge x_2)=Tx_1\wedge Tx_2$)  for all $x_1,x_2\in X$.  Vector lattices $X$
 and $Y$ are said to be \textit{lattice isomorphic\/} if there is a lattice isomorphism
 from $X$ onto $Y$, i.e., $T$ and $T^{-1}$ are lattice homomorphisms. Let $\Hom(X,Y)$ stand
 for the set of all lattice homomorphism from $X$ to $Y$. Recall also that the elements
 of the band $L_d^{\sim}(X,Y)\!:=\Hom(X,Y)^\perp$
 are referred to as \textit{diffuse operators}. An order bounded operator $T:X\to Y$ is said to be
 \textit{pseudoembedding\/} if $T$ belongs to the complementary band $L_a^{\sim}(X,Y)\!:=
 \Hom(X,Y)^{\perp\perp}$, the band generated by all lattice homomorphisms.
 Put $X^\sim\!:=L^\sim(X,\mathbb{R})$ and $X^\sim_a\!:=L^\sim_a(X,\mathbb{R})$.

 We~let $\!:=$ denote the~assignment by definition, while $\mathbb{N}$, $\mathbb{R}$,
 and $\mathbb{C}$ symbolize the~naturals, the~reals, and the~complexes. Throughout the sequel
 $\mathbb{B}$ is a complete Boolean algebra
 with top $\mathbb{1}$ and bottom $\mathbb{0}$.
 A~{\it partition of unity\/} in~$\mathbb{B}$ is a~family
 $(b_{\xi})_{\xi\in \Xi}$ in $\mathbb{B}$ such that
 $\bigvee_{\xi\in \Xi} b_{\xi}=\mathbb{1}$ and
 $b_{\xi}\wedge b_{\eta}=\mathbb{O}$ whenever $\xi\ne \eta$.

 The~reader can find the~relevant information on the~theory
 of order bounded operators in~\cite{AB, DOP, LZ, Vul,  Zaan}; on
 the Boolean valued models of set theory, in \cite{Bell, TZ}; and on Boolean valued analysis, in~\cite{BVA, IBA}.

 \section*{2. Pseudoembedding Operators}

 In this section we will give a~description of the band generated by disjointness
 preserving operators in the~vector lattice of order bounded operators.
 First we examine the scalar case.

 \subsec{2.1.}\proclaim{}For an arbitrary  vector lattice
 $X$ there exist a~unique cardinal $\gamma$ and a~disjoint family $(\varphi_ \alpha)_{ \alpha<\gamma}$ 
of nonzero 
 lattice homomor\-phisms $\varphi_ \alpha:X\rightarrow\mathbb{R}$
 such that every $f\in X^{\sim}$
 admits the unique representation
 $$
 f=f_d+\osum_{ \alpha<\gamma}\lambda_ \alpha\varphi_ \alpha
 $$
 where $f_d\in X_d^{\sim}$ and
 $(\lambda_ \alpha)_{ \alpha<\gamma}\subset\mathbb{R}$. The
 family $(\varphi_ \alpha)_{ \alpha<\gamma}$ is unique up to
 permutation and positive scalar multiplication.
 \Endproc

 \beginproof~The Dedekind complete vector lattice
 $X^{\sim}$ splits into the direct sum of the
 atomic band $X_a^{\sim}$ and the diffuse
 band $X_d^{\sim}\!:=(X_a^\sim)^\perp$; therefore, each functional
 $f\in E^{\sim}$ admits the unique representation $f=f_a+f_d$ with $f_a\in X_a^{\sim}$ and $f_d\in X_d^{\sim}$. Let $\gamma$ be the cardinality of the set $\mathcal{K}$ of one-dimensional bands
 in $X_a^{\sim}$ ($=$~atoms in
 $\mathbb{B}(X^{\sim})$).
 Then there exists a~family of lattice homomorphisms
 $(\varphi_ \alpha:X\rightarrow\mathbb{R})_{ \alpha<\gamma}$
 such that
 $\mathcal{K}=\{\varphi_ \alpha^{\perp\perp}:\, \alpha<\gamma\}$.
 It remains to observe that the mapping sending a~family of
 reals $(\lambda_ \alpha)_{ \alpha<\gamma}$ to the functional
 $x\mapsto\osum_{ \alpha<\gamma}\lambda_ \alpha\varphi_ \alpha(x)$
 implements a~lattice isomorphism between
 $X_a^{\sim}$ and some ideal in the vector
 lattice $\mathbb{R}^\gamma$.

 If $(\psi_ \alpha)_{ \alpha<\gamma}$ is a~disjoint family of nonzero 
 real lattice homomorphisms on $X$
 then for all $ \alpha,\beta<\gamma$ the functionals $\varphi_ \alpha$
 and $\psi_\beta$ are either disjoint or proportional with a~strictly positive coefficient,
 so that there exist a~permutation $(\omega_\beta)_{\beta<\gamma}$ of
 $(\varphi_ \alpha)_{ \alpha<\gamma}$ and a~unique family
 $(\mu_\beta)_{\beta<\gamma}$ in $\mathbb{R}_+$ such that
 $\psi_\beta=\mu_\beta\omega_\beta$ for all $\beta<\gamma$.~\endproof

 \subsec{2.2.}~Given two families $(S_ \alpha)_{ \alpha\in\mathrm{A}}$ and
 $(T_\beta)_{\beta\in\mathrm{B}}$ in $L^\sim(X,Y)$, say that
 $(S_ \alpha)_{ \alpha\in\mathrm{A}}$ is
 a~$\mathbb{P}(Y)$-\textit{permutation\/}
 of $(T_\beta)_{\beta\in\mathrm{B}}$ whenever there exists a~double family
 $(\pi_{ \alpha,\beta})_{ \alpha\in\mathrm{A},\,\beta\in\mathrm{B}}$
 in $\mathbb{P}(Y)$ such that $S_ \alpha=\sum_{\beta\in\mathrm{B}}\pi_{ \alpha,\beta}T_\beta$
 for all $ \alpha\in\mathrm{A}$, while
 $(\pi_{ \alpha,\bar{\beta}})_{ \alpha\in\mathrm{A}}$ and
 $(\pi_{\bar{ \alpha},\beta})_{\beta\in\mathrm{B}}$ are partitions
 of unity in $\mathbb{B}(Y)$ for all $\bar{ \alpha}\in\mathrm{A}$ and
 $\bar{\beta}\in\mathrm{B}$. It is easily seen that in case
 $Y=\mathbb{R}$ this amounts to saying that there is a~bijection
 $\nu:\mathrm{A}\rightarrow\mathrm{B}$ with
 $S_ \alpha=T_{\nu( \alpha)}$ for all $ \alpha\in\mathrm{A}$; i.e.,
 $(S_ \alpha)_{ \alpha\in\mathrm{A}}$ is a~permutation
 of $(T_\beta)_{\beta\in\mathrm{B}}$. We also say that
 $(S_ \alpha)_{ \alpha\in\mathrm{A}}$ is $\Orth(Y)$-multiple of
 $(T_ \alpha)_{ \alpha\in\mathrm{A}}$ whenever there exists a~family of orthomorphisms $(\pi_ \alpha)_{ \alpha\in\mathrm{A}}$ in $\Orth(Y)$ such that $S_ \alpha=\pi_ \alpha T_ \alpha$ for all $ \alpha\in\mathrm{A}$. In case $Y=\mathbb{R}$ we evidently get that $S_ \alpha$ is a~scalar multiple of $T_ \alpha$ for all $ \alpha\in\mathrm{A}$.

 Using the above notation, define the two mappings
 $\mathcal{S}:\mathrm{A}\rightarrow X^{{\scriptscriptstyle\wedge}\sim}{\downarrow}$ and
 $\mathcal{T}:\mathrm{B}\rightarrow X^{{\scriptscriptstyle\wedge}\sim}{\downarrow}$ by
 putting $\mathcal{S}( \alpha)\!:=S_ \alpha{\upwardarrow}$ $( \alpha\in\mathrm{A})$ and
 $\mathcal{T}(\beta)\!:=T_\beta{\upwardarrow}$
 $(\beta\in\mathrm{B})$. Recall that $\upwardarrow$ signify the modified ascent; see \cite[1.6.8]{BA_ST}.

 \subsec{2.3.}\proclaim{}Define the internal mappings $\tau,\sigma\in\mathbb{V}^{(\mathbb{B})}$ as
 $\sigma\!:=\mathcal{S}{\upwardarrow}$ and
 $\tau\!:=\mathcal{T}{\upwardarrow}$.
 Then $(\sigma( \alpha))_{ \alpha\in\mathrm{A}^{\scriptscriptstyle\wedge}}$
 is a~permutation of $(\tau(\beta))_{\beta\in\mathrm{B}^{\scriptscriptstyle\wedge}}$
 within~$\mathbb{V}^{(\mathbb{B})}$ if and only if $(S_ \alpha)_{ \alpha\in\mathrm{A}}$ is a~
 $\mathbb{P}(Y)$-\textit{permutation\/} of $(T_\beta)_{\beta\in\mathrm{B}}$. \Endproc

 \beginproof~Assume that $(\sigma( \alpha))_{ \alpha\in\mathrm{A}^{\scriptscriptstyle\wedge}}$
 is a~permutation of $(\tau(\beta))_{\beta\in\mathrm{B}^{\scriptscriptstyle\wedge}}$
 within $\mathbb{V}^{(\mathbb{B})}$.
 Then there is a~bijection $\nu:\mathrm{B}^{\scriptscriptstyle\wedge}
 \rightarrow\mathrm{A}^{\scriptscriptstyle\wedge}$ such that
 $\sigma( \alpha)=\tau(\nu( \alpha))$ for all
 $( \alpha\in\mathrm{A}^{\scriptscriptstyle\wedge})$.
 By \cite[1.5.8]{BA_ST} $\nu{\downwardarrow}$ is a~function from $\mathrm{A}$
 to $(\mathrm{B}^{\scriptscriptstyle\wedge}){\downarrow}=
 \mix(\{\beta^{\scriptscriptstyle\wedge}:\,\beta\in\mathrm{B}\})$.
 Thus, for each $ \alpha\in\mathrm{A}$ there exists a~partition of
 unity $(b_{ \alpha,\beta})_{\beta\in\mathrm{B}}$ such that
 $\nu{\downwardarrow}( \alpha)=\mix_{\beta\in\mathrm{B}}(b_{ \alpha,\beta}\beta^{\scriptscriptstyle\wedge})$. Since
 $\nu{\downwardarrow}$ is injective, we have
 $$
 \gathered
 \mathbb{1}=[\![(\forall \alpha_1, \alpha_2\in\mathrm{A}^{\scriptscriptstyle\wedge})
 (\nu( \alpha_1)=\nu( \alpha_2)\rightarrow \alpha_1= \alpha_2)]\!]
 \\
 =\bigwedge_{ \alpha_1, \alpha_2\in\mathrm{A}}
 [\![\nu( \alpha_1^{\scriptscriptstyle\wedge})
 =\nu( \alpha_2^{\scriptscriptstyle\wedge})\rightarrow
  \alpha_1^{\scriptscriptstyle\wedge}= \alpha_2^{\scriptscriptstyle\wedge}]\!]
 \\
 =\bigwedge_{ \alpha_1, \alpha_2}[\![\nu\downwardarrow( \alpha_1)
 =\nu\downwardarrow( \alpha_2)]\!]\Rightarrow
 [\![ \alpha_1^{\scriptscriptstyle\wedge}= \alpha_2^{\scriptscriptstyle\wedge}]\!],
 \endgathered
 $$
 and so $[\![\nu\downwardarrow( \alpha_1)=\nu\downwardarrow( \alpha_2)]\!]
 \leq[\![ \alpha_1^{\scriptscriptstyle\wedge}= \alpha_2^{\scriptscriptstyle\wedge}]\!]$
 for all $ \alpha_1, \alpha_2\in\mathrm{A}$. Taking this inequality and the definition
 of  $\nu{\downwardarrow}$ into account yields
 $$
 \gathered
 b_{ \alpha_1,\beta}\wedge b_{ \alpha_2,\beta}\leq[\![\nu{\downwardarrow}( \alpha_1)=\beta^{\scriptscriptstyle\wedge}]\!]
 \wedge[\![\nu{\downwardarrow}( \alpha_2)=\beta^{\scriptscriptstyle\wedge}]\!]
 \\
 \leq[\![\nu\downwardarrow( \alpha_1)=\nu\downwardarrow( \alpha_2)]\!]
 \leq[\![ \alpha_1^{\scriptscriptstyle\wedge}= \alpha_2^{\scriptscriptstyle\wedge}]\!],
 \endgathered
 $$
 so that $ \alpha_1\ne \alpha_2$ implies $b_{ \alpha_1,\beta}\wedge b_{ \alpha_2,\beta}=\mathbb{0}$ (because
 $x \ne y\Longleftrightarrow[\![x^{\scriptscriptstyle\wedge}=y^{\scriptscriptstyle\wedge}]\!]=\mathbb{O}$
 by \cite[1.4.5\,(2)]{BA_ST}. At the same time, surjectivity of $\nu$ implies
 $$
 \gathered
 \mathbb{1}=[\![(\forall\,\beta\in\mathrm{B}^{\scriptscriptstyle\wedge})
 (\exists\, \alpha\in\mathrm{A}^{\scriptscriptstyle\wedge})\beta=\nu( \alpha)]\!]
 \\
 =\bigwedge_{\beta\in\mathrm{B}}\bigvee_{ \alpha\in\mathrm{A}}[\![\beta^{\scriptscriptstyle\wedge}
 =\nu{\downwardarrow}( \alpha)]\!]=\bigwedge_{\beta\in\mathrm{B}}\bigvee_{ \alpha\in\mathrm{A}}
 b_{ \alpha,\beta}.
 \endgathered
 $$
 It follows that $(b_{ \alpha,\beta})_{ \alpha\in\mathrm{A}}$ is a~partition of unity in
 $\mathbb{B}$ for all $\beta\in\mathrm{B}$. By the choice of $\nu$ it follows that
 $b_{ \alpha,\beta}\leq[\![\sigma( \alpha^{\scriptscriptstyle\wedge})=\tau(\beta^{\scriptscriptstyle\wedge})]\!]$,
 because of the estimations
 $$
 \gathered
 b_{ \alpha,\beta}\leq[\![\sigma( \alpha^{\scriptscriptstyle\wedge})=
 \tau(\nu( \alpha^{\scriptscriptstyle\wedge}))]\!]\wedge[\![\nu( \alpha^{\scriptscriptstyle\wedge})=\beta^{\scriptscriptstyle\wedge}]\!]
 \\
 \leq[\![\sigma( \alpha^{\scriptscriptstyle\wedge})=
 \tau(\beta^{\scriptscriptstyle\wedge})]\!]
 =[\![\mathcal{S}( \alpha)=\mathcal{T}(\beta)]\!].
 \endgathered
 $$
 Put now $\pi_{ \alpha,\beta}\!:=\chi(b_{ \alpha,\beta})$ and observe that
 $b_{ \alpha,\beta}\leq[\![\mathcal{S}( \alpha)x^{\scriptscriptstyle\wedge}=
 \mathcal{T}(\beta)x^{\scriptscriptstyle\wedge}]\!]
 \leq[\![S_ \alpha x=T_\beta x]\!]$ for all  $ \alpha\in\mathrm{A}$, $\beta\in\mathrm{B}$, and $x\in X$.
 Using $(\mathbb{G})$, we obtain $\pi_{ \alpha,\beta}S_ \alpha=\pi_{ \alpha,\beta}T_\beta$ and so
 $S_ \alpha=\sum_{\beta\in\mathrm{B}}\pi_{ \alpha,\beta}T_\beta$ for all $ \alpha\in\mathrm{A}$. Clearly,
 $(\pi_{ \alpha,\beta})$ is the family as required in Definition
 2.2. The sufficiency is shown by the same reasoning in the reverse direction.~\endproof

 \subsec{2.4.}~A nonempty set $\mathcal{D}$ of positive operators from $X$ to $Y$ is called {\it strongly
 generating\/} if $\mathcal{D}$ is a disjoint set and $S(X)^{\perp\perp}=Y$ for all $S\in\mathcal{D}$. If, in addition,
 $\mathcal{D}^{\perp\perp}=B$, then we say also that $\mathcal{D}$ \textit{strongly
 generates\/} the band $B\subset L^{\sim}(X,Y)$ or $B$ is strongly generated by $\mathcal{D}$.
 In case $Y=\mathbb{R}$, the strongly generating sets in $X^\sim=L^{\sim}(X,\mathbb{R})$
 are precisely disjoint sets of nonzero  positive functionals.

 Given a~cardinal $\gamma$ and a~universally complete vector lattice $Y$, say that
 a~vector lattice $X$ is $(\gamma,Y)$-\textit{homogeneous\/} if the band
 $L^{\sim}_a(X,Y)$ is strongly generated by a~set of lattice
 homomorphisms of cardinality $\gamma$ and for every nonzero projection $\pi\in\mathbb{P}(Y)$
 and every strongly generating set $\mathcal{D}$ in $L^{\sim}_a(X,\pi Y)$ we have $\card(\mathcal{D})\geq\gamma$.
 We say also that $X$ is $(\gamma,\pi)$-homogeneous if $\pi\in\mathbb{P}(Y)$ and $X$
 is $(\gamma,\pi Y)$-homogeneous.
 Evidently, the $(\gamma,\mathbb{R})$-homogeneity of a~vector lattice $X$  amounts just
 to saying that the band $X_a^{\sim}$ is generated in $X^\sim$ by
 a~cardinality $\gamma$ disjoint set of nonzero   lattice homomorphisms  or, equivalently,
 the cardinality of the set of
 atoms in $\mathbb{B}(X^\sim)$ equals $\gamma$.

 Take $\mathcal{D}\subset L^\sim(X,\mathcal{R}{\downarrow})$ and  $\Delta\in\mathbb{V}^{(\mathbb{B})}$ with $[\![\Delta\subset (X^{\scriptscriptstyle\wedge})^\sim]\!]=\mathbb{1}$. Put
 $\mathcal{D}_{\uparrow}\!:=\{T{\upwardarrow}:\,T\in\mathcal{D}\}{\uparrow}$ and $\Delta^{\downarrow}\!:=\{\tau{\downwardarrow}:\,\tau\in\Delta{\downarrow}\}$.
 Let $\mix(\mathcal{D})$ stand for the set of all $T\in L^\sim(X,\mathcal{R}{\downarrow})$ representable as $Tx=\osum_{\xi\in\Xi}\pi_\xi T_\xi x$ $(x\in X)$ with $(\pi_\xi)_{\xi\in\Xi}$ a~partition of unity  in $\mathbb{P}(\mathcal{R}{\downarrow})$ and  $(T_\xi)_{\xi\in\Xi}$ a~family in $\mathcal{D}$.

 \subsec{2.5}~\proclaim{}Let $\Delta\subset(X^{\scriptscriptstyle\wedge})^\sim$ be 
a~disjoint set of nonzero  positive func\-tio\-nals which has cardinality $\gamma^{\scriptscriptstyle\wedge}$ within $\mathbb{V}^{(\mathbb{B})}$. Then there exists a~cardinality $\gamma$ strongly generating set of positive operators $\mathcal{D}$ from $X$ to $\mathcal{R}{\downarrow}$ such that
 $\Delta=\mathcal{D}_\uparrow$ and $\Delta^\downarrow=\mix(\mathcal{D})$.
 \Endproc

 \beginproof~If $\Delta$ obeys the conditions then there is $\phi\in\mathbb{V}^{(\mathbb{B})}$ such that $[\![\phi:\gamma^{\scriptscriptstyle\wedge}\rightarrow\Delta$ is a~bijection$]\!]=\mathbb{1}$. Note that $\varphi{\downwardarrow}$ sends $\gamma$ into $\Delta{\downarrow}\subset(X^{\scriptscriptstyle\wedge})^\sim{\downarrow}$
 by \cite[1.5.8]{BA_ST}. By \cite[Theorem 3.3.3]{BA_ST}, we can define the mapping $ \alpha\mapsto\Phi( \alpha)$ from $\gamma$ to $L^\sim(X,\mathcal{R}{\downarrow})$ by putting $\Phi( \alpha)\!:=(\phi{\downwardarrow}( \alpha)){\downwardarrow}$.
 Put $\mathcal{D}\!:=\{\Phi( \alpha):\, \alpha\in\gamma\}$ and note that $\mathcal{D}\subset\Delta^{\downarrow}$.
 Using~\cite[1.6.6]{BA_ST} and the surjectivity of $\phi$ we have $\Delta{\downarrow}=
 \varphi(\gamma^{\scriptscriptstyle\wedge}){\downarrow}=
 \mix\{\phi{\downwardarrow}( \alpha)):\, \alpha\in\gamma\}$ and combining this with  \cite[3.3.7]{BA_ST} we get
 $\Delta=\mathcal{D}_\uparrow$ and $\Delta^\downarrow=\mix(\mathcal{D})$.

 The injectivity of $\phi$ implies that
 $$
 [\![(\forall \alpha,\beta\in\gamma^{\scriptscriptstyle\wedge})( \alpha\ne\beta
 \rightarrow\phi( \alpha)\ne\phi(\beta)]\!]=\mathbb{1}.
 $$
  Replacing the universal quantifier by
 the supremum over $ \alpha,\beta\in\gamma^{\scriptscriptstyle\wedge}$, from \cite[1.4.5\,(1) and 1.4.5\,(2)]{BA_ST} we deduce that
 $$
 \mathbb{1}=\bigwedge_{ \alpha,\beta\in\gamma}
 [\![ \alpha^{\scriptscriptstyle\wedge}\ne\beta^{\scriptscriptstyle\wedge}]\!]
 \Rightarrow[\![\varphi( \alpha^{\scriptscriptstyle\wedge})\ne\phi(\beta)^{\scriptscriptstyle\wedge}]\!]
 =\bigwedge_{\substack{{ \alpha,\beta\in\gamma}\\{ \alpha\ne\beta}}}
 [\![\Phi( \alpha)\ne\Phi(\beta)]\!],
 $$
 and so  $ \alpha\ne\beta$ implies $\Phi( \alpha)\ne\Phi(\beta)$ for all $ \alpha,\beta\in\gamma$.
 Thus $\Phi$ is injective and the cardinality of $\mathcal{D}$ is $\gamma$.
 The fact that $\mathcal{D}$ is strongly generating follows from \cite[3.3.5\,(5) and  3.8.4]{BA_ST}.~\endproof

 \subsec{2.6.}\proclaim{}If $\mathcal{D}$ is a~cardinality $\gamma$ strongly generating set of positive operators
 from $X$ to $\mathcal{R}{\downarrow}$ of  then $\Delta=\mathcal{D}_\uparrow\subset
 (X^{\scriptscriptstyle\wedge})^\sim$ is a~disjoint set of nonzero 
 positive func\-tio\-nals which has
 cardinality $|\gamma^{\scriptscriptstyle\wedge}|$ within $\mathbb{V}^{(\mathbb{B})}$.
 \Endproc

 \beginproof~Assume that $\mathcal{D}\subset L(X,\mathcal{R}{\downarrow})$ is a~strongly generating
 set of car\-di\-na\-lity $\gamma$. Then there is a~bijection $f:\gamma\rightarrow\mathcal{D}{\upwardarrow}$.
 Moreover, $ \alpha\ne\beta$ implies $[\![f( \alpha)\perp f(\beta)]\!]=\mathbb{1}$ by \cite[3.3.5\,(5)]{BA_ST}
 and $[\![f( \alpha)\ne0]\!]=\mathbb{1}$ by~\cite[3.8.4]{BA_ST}. Interpreting in $\mathbb{V}^{(\mathbb{B})}$ the $\ZFC$-theorem
 $$
 (\forall f,g\in X^\sim)(f\ne0\wedge g\ne0\wedge f\perp g\rightarrow f\ne g)
 $$
 yields $[\![f( \alpha)\ne f(\beta)]\!]=\mathbb{1}$ for all $ \alpha,\beta\in\gamma$, $ \alpha\ne\beta$.
 It follows that $\phi\!:=f{\upwardarrow}$ is a~bijection from $\gamma^{\scriptscriptstyle\wedge}$
 onto $\Delta=(\mathcal{D}{\upwardarrow}){\uparrow}$, so that the cardinality of $\Delta$ is
 $|\gamma^{\scriptscriptstyle\wedge}|$. The proof is completed by the arguments similar to
 those in~2.5.~\endproof

 \subsec{2.7.}\proclaim{}A vector lattice $X$ is
 $(\gamma,\mathcal{R}{\downarrow})$-homogeneous for some cardinal $\gamma$ if and only if
 $[\![\,\gamma^{\scriptscriptstyle\wedge}$ is a~cardinal and
 $X^{\scriptscriptstyle\wedge}$ is
 $(\gamma^{\scriptscriptstyle\wedge},\mathcal{R})$-homogeneous
 $]\!]=\mathbb{1}$.
 \Endproc

 \beginproof~\textit{Sufficiency:}~Assume that  $\gamma^{\scriptscriptstyle\wedge}$ is a~cardinal and
 $X^{\scriptscriptstyle\wedge}$ is $(\gamma^{\scriptscriptstyle\wedge},\mathcal{R})$-homo\-ge\-neous
 within $\mathbb{V}^{(\mathbb{B})}$. The latter means that
 $(X^{\scriptscriptstyle\wedge})_a^\sim$ is generated by a~disjoint set of nonzero 
 lattice
 homomorphisms  $\Delta\subset(X^{\scriptscriptstyle\wedge})^\sim$  which has cardinality
 $\gamma^{\scriptscriptstyle\wedge}$ within $\mathbb{V}^{(\mathbb{B})}$. By~2.5
 there exists a~strongly generating set $\mathcal{D}$ in $L^\sim_a(X,\mathcal{R}{\downarrow})$
 of cardinality $\gamma^{\scriptscriptstyle\wedge}$ such that $\Delta=\mathcal{D}_{\uparrow}$.
 Take a~nonzero $\pi\in\mathbb{P}(\mathcal{R}{\downarrow})$ and put $b\!:=\chi^{-1}(\pi)$. Recall
 that  we can identify $L^\sim(X,\pi(\mathcal{R}{\downarrow}))$ and $L^\sim(X,(b\wedge\mathcal{R}){\downarrow})$.
 If $\mathcal{D}'$ is a~strongly generating set in $L_a^\sim(X,\pi(\mathcal{R}{\downarrow}))$
 of cardinality $\beta$ then $\mathcal{D}^\prime_{\uparrow}$ strongly generates
 $(X^{\scriptscriptstyle\wedge})_a^\sim$ and has  cardinality $|\beta^{\scriptscriptstyle\wedge}|$
 within the relative universe $\mathbb{V}^{([\mathbb{0},b])}$. By \cite[1.3.7]{BA_ST}
 $\gamma^{\scriptscriptstyle\wedge}=
 |\beta^{\scriptscriptstyle\wedge}|\leq \beta^{\scriptscriptstyle\wedge}$ and so $\gamma\leq\beta$.

 \textit{Necessity:}~Assume now that $X$ is $(\gamma,\mathcal{R}{\downarrow})$-homogeneous and
 the~set $\mathcal{D}$ of lattice homomorphisms  of cardinality $\gamma$ generates strongly the
 band $L_a^\sim(X,\mathcal{R}{\downarrow})$. Then $\Delta=\mathcal{D}_{\uparrow}$ generates
 the band $(X^{\scriptscriptstyle\wedge})_a^\sim$ and the cardinalities  of $\Delta$ and
 $\gamma^{\scriptscriptstyle\wedge}$ coincide; i.e., $|\Delta|=|\gamma^{\scriptscriptstyle\wedge}|$.
 By \cite[1.9.11]{BA_ST} the cardinal  $|\gamma^{\scriptscriptstyle\wedge}|$ has the representation
 $|\gamma^{\scriptscriptstyle\wedge}| =\mix_{ \alpha\leq\gamma}b_ \alpha \alpha^{\scriptscriptstyle\wedge}$,
 where $(b_ \alpha)_{ \alpha\leq\gamma}$ is a~partition of unity in $\mathbb{B}$.
 It follows that $b_ \alpha\leq[\![\Delta$ is a~generating set in $(X^{\scriptscriptstyle\wedge})_a^\sim$
 of cardinality $ \alpha^{\scriptscriptstyle\wedge}]\!]=\mathbb{1}$.
 If $b_ \alpha\ne\mathbb{0}$ then $b_ \alpha\wedge\Delta$ is a~generating
 set in $(X^{\scriptscriptstyle\wedge})_a^\sim$ of cardinality $|\gamma^{\scriptscriptstyle\wedge}|=
  \alpha^{\scriptscriptstyle\wedge}\leq\gamma^{\scriptscriptstyle\wedge}$
  in the relative universe $\mathbb{V}^{[\mathbb{0},b_ \alpha]}$.
 Put $\pi_ \alpha=\chi(b_ \alpha)$ and $\pi_ \alpha\circ\mathcal{D}\!:=\{\pi_ \alpha\circ T:\,T\in\mathcal{D}\}$.
 Clearly, $b_ \alpha\wedge\Delta=(\pi_ \alpha\mathcal{D})_{\uparrow}$ and so $\pi_ \alpha\circ\mathcal{D}$ strongly generates the band $L_a^\sim(X,\mathcal{R}{\downarrow})$. By hypothesis $\mathcal{D}$ is $(\gamma,\mathcal{R}{\downarrow})$-homogeneous, consequently, $ \alpha\geq\gamma$, so that $ \alpha=\gamma$, since $ \alpha\leq\gamma$ if and only if $ \alpha^{\scriptscriptstyle\wedge}\leq\gamma^{\scriptscriptstyle\wedge}$.
 Thus, $|\gamma^{\scriptscriptstyle\wedge}|=\gamma^{\scriptscriptstyle\wedge}$
 whenever $b_ \alpha\ne\mathbb{0}$ and $\gamma^{\scriptscriptstyle\wedge}$ is a~cardinal within $\mathbb{V}^{(\mathbb{B})}$.~\endproof

 \subsec{2.8.}\proclaim{}Let $X$ be a~$(\gamma,Y)$-homogeneous vector lattice for some
 universally complete vector lattice $Y$ and a~nonzero cardinal $\gamma$. Then there exists
 a~strongly generating family of lattice homomorphisms $(\Phi_{ \alpha})_{ \alpha<\gamma}$
 from $X$ to $Y$ such that each operator $T\in L_a^{\sim}(X,Y)$ admits the unique representation
 $T=\osum_{ \alpha<\gamma}\sigma_{ \alpha}\circ\Phi_{\gamma, \alpha}$,
 where $(\sigma_{ \alpha})_{ \alpha<\gamma}$ is a~family of orthomorphisms in $\Orth(Y)$.
 \Endproc

 \beginproof~This is immediate from the definitions in 2.4.~\endproof

 \subsec{2.9.}~\theorem{}Let $X$ and $Y$ be vector lattices with $Y$ universally complete. Then there are
 a~nonempty set of cardinals $\Gamma$ and a~partition of unity $(Y_\gamma)_{\gamma\in\Gamma}$
 in $\mathbb{B}(Y)$ such that $X$ is $(\gamma,Y_\gamma)$-homogeneous for all $\gamma\in\Gamma$.
 \Endproc

 \beginproof~We may assume without loss of generality that
 $Y=\mathcal{R}{\downarrow}$. The transfer principle tells us
 that according to~2.1 there exists a~cardinal $\varkappa$ within $\mathbb{V}^{(\mathbb{B})}$ such that
 $(X^{\scriptscriptstyle\wedge})^\sim_a$ is generated by a~cardinality~$\varkappa$
  disjoint set $\mathcal{H}$ of nonzero  $\mathbb{R}^{\scriptscriptstyle\wedge}$-linear lattice homomorphisms
  or, equivalently,
 $[\![X^{\scriptscriptstyle\wedge}\text{~is~}(\varkappa, \mathcal{R})\text{-homogeneous}\,]\!]=\mathbb{1}$.
 By \cite[1.9.11]{BA_ST} there is a~nonempty set of cardinals $\Gamma$ and a~partition of unity $(b_\gamma)_{\gamma\in\Gamma}$ in $\mathbb{B}$ such that $\varkappa=\mix_{\gamma\in\Gamma}
 b_\gamma\gamma^{\scriptscriptstyle\wedge}$.
 It follows that $b_\gamma\leq[\![X^{\scriptscriptstyle\wedge}
 \text{~is~}(\gamma^{\scriptscriptstyle\wedge},\mathcal{R})
 \text{-homogeneous}\,]\!]$ for all $\gamma\in\Gamma$. Passing to the relative subalgebra
 $\mathbb{B}_\gamma\!:=[\mathbb{0},b_\gamma]$ and considering \cite[1.3.7]{BA_ST} we conclude that
 $\mathbb{V}^{(\mathbb{B}_\gamma)}\models ``X^{\scriptscriptstyle\wedge}\text{~is~}
 (\gamma^{\scriptscriptstyle\wedge},b_\gamma\wedge\mathcal{R})
 \text{-homogeneous''}$, so that $X$ is $(\gamma,(b_\gamma\wedge\mathcal{R}){\downarrow})$-homogeneous by 2.7.
 In view of \cite[2.3.6]{BA_ST} $(b_\gamma\wedge\mathcal{R}){\downarrow}$ is lattice isomorphic to $Y_\gamma$,
 and so the desired result follows.~\endproof

 \subsec{2.10.}~\theorem{}Let $X$ and $Y$ be vector lattices
 with $Y$ uni\-ver\-sally complete. Then there are a~nonempty set of cardinals  $\Gamma$ and
 a~partition of unity $(Y_\gamma)_{\gamma\in\Gamma}$ in $\mathbb{B}(Y)$ such that
 to each cardinal $\gamma\in\Gamma$ there is a~disjoint family
 of  lattice homomorphisms $(\Phi_{\gamma, \alpha})_{ \alpha<\gamma}$ from $X$ to $Y_\gamma$ satisfying
 \vspace*{2pt}

 \subsec{\bf(1)}~$\Phi_{\gamma, \alpha}(X)^{\perp\perp}=Y_\gamma\ne\{0\}$ for all $\gamma\in\Gamma$ and $ \alpha<\gamma$.
 \vspace*{2pt}

 \subsec{\bf(2)}~$X$ is $(\gamma,Y_\gamma)$-homogeneous for all $\gamma\in\Gamma$.
 \vspace*{2pt}

 \subsec{\bf(3)}~For each order dense sublattice $Y_0\subset Y$ each
 $T\in L^{\sim}(X,Y_0)$ admits the unique
 representation
 $$
 T=T_d+\osum_{\gamma\in\Gamma}
 \osum_{ \alpha<\gamma}\sigma_{\gamma, \alpha}\circ\Phi_{\gamma, \alpha},
 $$
 with $T_d\in L^{\sim}_d(X,Y)$ and
 $\sigma_{\gamma, \alpha}\in\Orth(\Phi_{\gamma, \alpha},
 Y_0)$.

 For every $\gamma\in\Gamma$ the family $(\Phi_{\gamma, \alpha})_{ \alpha<\gamma}$
 is unique up to $\mathbb{P}(Y)$-permutation and $\Orth(Y_{\gamma})_+$-multiplication.
 \Endproc

 \beginproof~The existence of $(Y_\gamma)_{\gamma\in\Gamma}$ and $(\Phi_{\gamma,
 \alpha})_{\gamma\in\Gamma,\, \alpha<\gamma}$ with the required properties is immediate
 from 2.8 and 2.9. The uniqueness follows from 2.1 and 2.3.~\endproof

 \subsec{2.11.}~Theorem 2.10, the main result of Section 2, was proved by Tabuev in~\cite[Theorem 2.2]{Tab2} with standard tools. The pseudoembedding operators are closely
 connected with the so-called order narrow operators. A linear operator $T:X\to Y$ is
 \textit{order narrow\/} if for every $x\in X_+$ there exists a~net $(x_ \alpha)$ in $X$
 such that $|x_ \alpha|=x$ for all $ \alpha$ and $(Tx_ \alpha)$ is order convergent to zero
 in $Y$; see \cite[Definition 3.1]{MMP}. The main result by
 Maslyuchenko, Mykhaylyuk, and Popov in~\cite[Theorem 11.7\,(ii)]{MMP} states that if $X$ and $Y$ are Dedekind complete
 vector lattices with $X$ atomless and $Y$  an order ideal of some order continuous
 Banach lattice then an order bounded order continuous operator is order narrow if and only
 if it is pseudoembedding.

 \subsec{2.12.}~The term {\it pseudoembedding operator} stems from a~result by
 Rosenthal \cite{Ros} which asserts that a~nonzero bounded linear operator in $L^1$ is a~pseudoembedding if and only
 if it is a~near isometric embedding when restricted to a~suitable $L^1(A)$-subspace. Systematic study of
 narrow operators was started by Plichko and Popov in \cite{PP}. Concerning a~detailed presentation of the theory of narrow operators, we refer to the recent book by Popov and Randrianantoanina \cite{PR} and the references therein.

 \section*{3.~The Noncommutative Wickstead Problem}

 \textsl{When are we so happy in a~vector lattice that all band preserving linear operators turn out to be order bounded?}
  
  This question was raised by~Wickstead in~\cite{Wic1}. We refer to \cite{IBA, BA_ST} and \cite{GKK} for a detailed presentation of the Wickstead problem. In this section we consider a noncommutative version of the problem. The relevant information on the theory of Baer $\ast$-algebras and $AW^ \ast$-algebras  can be found
 in~Berberian~\cite{Berb}, Chilin~\cite{Chil}, Kusraev~\cite{DOP}.

 \subsec{3.1.}~A {\it Baer $ \ast$-algebra\/} is a~complex
 involutive algebra $A$ such that, for each nonempty
 $M\subset A$, there is a~projection, i.e., a~hermitian
 idempotent $p$, satisfying $M^\perp=pA$ where
 $M^\perp\!:=\{y\in A:(\forall\,x\in M)\,xy=0\}$ is the
 right annihilator of~$M$. Clearly,
 this amounts to saying that each left annihilator
 has the form ${}^\perp M=Aq$ for an appropriate projection
 $q$. To each left annihilator $L$ in a~Baer \mbox{$ \ast$-algebra}
 there is a~unique projection $q_L\in A$ such that $x=xq_L$
 for all $x\in L$ and $q_L y=0$ whenever
 $y\in L^\perp$. The mapping $L\mapsto q_L$ is an isomorphism
 between the poset of left annihilators and the poset of
 all projections. Thus, the poset $\mathbb{P}(A)$ of all projections in a~Baer
 \mbox{$ \ast$-algebra} is an~order complete lattice. \big(Clearly, the formula
 $q\leq p\Longleftrightarrow q=qp=pq$, sometimes pronounced  as  ``$p$ contains
 $q$,'' specifies some order on the set of
 projections $\mathbb{P}(A)$.\big)

  An element $z$ in $A$ is  {\it central\/} provided that $z$ commutes with every member of~$A$; i.e.,
 $(\forall\,x\in A)\,xz=zx$. The {\it center\/} of a~Baer $ \ast$-algebra  $A$ is the set $\mathcal{Z}(A)$
 comprising central elements. Clearly, $\mathcal{Z}(A)$
 is a~commutative Baer $ \ast$-subalgebra of~$A$, with
 $\lambda\mathbb{1}\in\mathcal{Z}(A)$ for all
 $\lambda\in\mathbb{C}$. A {\it central projection\/}
 of $A$ is a~projection belonging to $\mathcal{Z}(A)$. Put
 $\mathbb{P}_c(A)\!:=\mathbb{P}(A)\cap\mathcal{Z}(A)$.

 \subsec{3.2.}~A {\it derivation\/} on a~Baer $ \ast$-algebra $A$
 is a~linear operator $d:A\to A$ satisfying $d(xy)=d(x)y+xd(y)$
 for all $x,y\in A$. A derivation $d$ is  {\it
 inner\/} provided that $d(x)=ax-xa$ $(x\in A)$ for some $a\in A$. Clearly, an
 inner derivation vanishes on $\mathcal{Z}(A)$ and is
 $\mathcal{Z}(A)$-linear; i.e., $d(ex)=ed(x)$ for all $x\in A$
 and $e\in\mathcal{Z}(A)$.

 Consider a~derivation $d:A\to A$ on a~Baer
 $ \ast$-algebra $A$. If $p\in A$ is a~central projection
 then $d(p)=d(p^2)=2pd(p)$. Multiplying this identity
 by $p$ we have $pd(p)=2pd(p)$ so that $d(p)=pd(p)=0$.
 Consequently, every derivation vanishes on the linear span
 of $\mathbb{P}_c(A)$, the set of all central projections.
 In particular, $d(ex)=ed(x)$
 whenever $x\in A$ and $e$ is a~linear combination of central
 projections. Even if the linear span of central projections
 is dense in a~sense  in $\mathcal{Z}(A)$, the derivation $d$
 may fail to be $\mathcal{Z}(A)$-linear.

 This brings up the natural question: \textsl{Under what
 conditions is every derivation
 $Z$-linear on a~Baer $ \ast$-algebra $A$  provided that $Z$ is a~Baer $ \ast$-subalgebra
 of $\mathcal{Z}(A)$?}

 \subsec{3.3.}~An~{\it $AW^ \ast$-algebra\/} is a~$C^ \ast$-algebra with unity $\mathbb{1}$ which is also a~Baer
 $ \ast$\hbox{-}algebra. More explicitly, an~$AW^ \ast$-algebra is
 a~$C^ \ast$-algebra whose every right annihilator has the
 form~$pA$, with $p$ a~projection. Clearly, $\mathcal{Z}(A)$
 is a~commutative $AW^*$\hbox{-}subalgebra of~$A$. If $\mathcal{Z}(A)=
 \{\lambda\mathbb{1}:\lambda\in\mathbb{C}\}$ then the
 $AW^*$-algebra $A$ is  an~$AW^*$\hbox{-}{\it factor}.

 \subsec{3.4.}\Proclaim{}A~$C^ \ast$-algebra $A$ is an~$AW^ \ast$-algebra
 if and only if the two conditions hold:

 \subsec{(1)} Each orthogonal family in $\mathbb{P}(A)$
 has a~supremum.

 \subsec{(2)} Each maximal commutative $ \ast$-subalgebra
 of $A_0\subset A$ is a~Dedekind complete $f$-algebra {\rm(}or,
 equivalently, coincides with the least norm closed
 $ \ast$-subalgebra containing all  projections of~$A_0)$.
 \Endproc

 \subsec{3.5.}~Given an $AW^ \ast$-algebra $A$,
 define the two sets $C(A)$ and $S(A)$ of measurable and locally
 measurable operators, respectively. Both are Baer
 $ \ast$-algebras; cp. Chilin~\cite{Chil}. Suppose that $\Lambda$ is an
 $AW^\ast$-sub\-al\-gebra in $\mathcal{Z}(A)$, and $\Phi$ is
 a~$\Lambda$-valued trace on $A_+$. Then we can define another
 Baer $ \ast$-algebra, $L(A,\Phi)$, of $\Phi$-measurable
 operators. The center $\mathcal{Z}(A)$ is a~vector lattice
 with a~strong unit, while the centers of
 $C(A)$, $S(A)$, and $L(A,\Phi)$ coincide with the universal
 completion of~$\mathcal{Z}(A)$.
 If $d$ is a~derivation  on $C(A)$, $S(A)$, or
 $L(A,\Phi)$ then $d(px)=pd(x)$ $\bigl(p\in\mathbb{P}_c(A)\bigr)$ so
 that $d$ can be considered as band preserving in a~sense
 (cp. \cite[4.1.1 and 4.10.4]{BA_ST}). The natural question arises immediately
 about these  algebras:

 \subsec{3.6.}~\textrm{WP(C)}:~\textsl{When are all derivations
 on $C(A)$, $S(A)$, or $L(A,\Phi)$  inner?} This question may be
 regarded as the \textit{noncommutative Wickstead
 problem}.

 \subsec{3.7.}~The classification of $AW^ \ast$-algebras
 into types is determined from  their lattices
 of projections $\mathbb{P}(A)$; see Sakai~\cite{Sak}. We recall  only the definition
 of type~I $AW^ \ast$-algebra. A~projection $\pi\in A$ is  {\it abelian\/} if  $\pi A\pi$ is a~commutative algebra. An algebra $A$ has {\it type\/}~I provided that each  nonzero projection in $A$ contains a~nonzero abelian  projection.

 A~$C^ \ast$-algebra $A$ is  ${\mathbb B}$-{\it embeddable}
 provided that there are a~type~I $AW^ \ast$-algebra $N$ and a~$ \ast$-monomorphism
 $\imath: A\rightarrow N$ such that ${\mathbb B}={\mathbb P}_c(N)$
 and $\imath (A)=\imath(A)^{\prime\prime}$, where
 $\imath (A)^{\prime\prime}$ is the bicommutant
 of~$\imath (A)$ in~$N$. Note that in this event $A$ is
 an~$AW^ \ast$-algebra and ${\mathbb B}$ is a~complete
 subalgebra of~${\mathbb P}_c(A)$.

 \subsec{3.8.} \theorem{}Let $A$ be a~type I
 $AW^ \ast$-algebra, let $\Lambda$ be an $AW^ \ast$\hbox{-}subalgebra of
 $\mathcal{Z}(A)$, and let $\Phi$ be a~$\Lambda$-valued faithful
 normal semifinite trace on $A$. If the complete Boolean algebra
 $\mathbb{B} \!:=\mathbb{P}(\Lambda)$ is $\sigma$-distributive
 and $A$ is $\mathbb{B}$-embeddable, then every derivation
 on $L(A,\Phi)$ is inner.
 \Endproc

 \beginproof~We briefly sketch the proof.
 Let $\mathcal{A}\in\mathbb{V}^{(\mathbb{B})}$
 be the Boolean valued representation of~$A$. Then $\mathcal{A}$
 is a~von Neumann algebra within $\mathbb{V}^{(\mathbb{B})}$.
 Since  the Boolean valued interpretation preserves classification
 into types, $\mathcal{A}$ is of type I. Let $\varphi$ stand for
 the Boolean valued representation of~$\Phi$. Then $\varphi$
 is a faithful normal semifinite $\mathcal{C}$-valued
 trace on $\mathcal{A}$ and the descent of $L(\mathcal{A},\varphi)$
 is $ \ast$-$\Lambda$-isomorphic to $L(A,\Phi)$; cp.~Korol$'$ and Chilin~\cite{KCh}.
 Suppose that $d$ is a~derivation on $L(A,\Phi)$ and
 $\delta$ is the Boolean valued representation~of~$d$. Then $\delta$
 is a~$\mathcal{C}$-valued $\mathbb{C}^{\scriptscriptstyle\wedge}$-linear
 derivation on $L(\mathcal{A},\varphi)$. Since $\mathbb{B}$ is
 $\sigma$-distributive, $\mathcal{C}=
 {\mathbb C}^{\scriptscriptstyle\wedge}$ within
 $\mathbb{V}^{(\mathbb{B})}$ and $\delta$ is
 $\mathcal{C}$-linear. But it is well known that every derivation
 on a~type I von Neumann algebra is inner; cp. \cite{AAKud}.
 Therefore, $d$ is also inner.~\endproof

 \subsec{3.8.}~Theorem 3.8 is taken from Gutman, Kusraev, and Kutateladze~\cite[Theorem 4.3.6]{GKK}. This fact
  is an interesting ingredient of the theory of noncommutative integration which stems from Segal \cite{Seg}. Considerable attention is given to derivations on various algebras of measurable operators associated with an $AW^ \ast$-algebra and a~central-valued trace. We mention only the article \cite{AAKud} by Albeverio, Ajupov, and Kudaybergenov and the article \cite{BPS}
 by Ber, de~Pagter, and Sukochev.

 \section*{4.~The Radon--Nikod\'ym Theorem for \emph{JB}-Algebras}

  In this section we sketch some further applications of the Boolean value approach to a nonassociative Radon--Nikod\'ym type  theorems.

 \subsec{4.1.}~Let $A$ be a~vector space over some field $\mathbb{F}$. Say that $A$
 is a~{\it Jordan algebra,} if there is given a (generally) nonassociative
 binary operation $A\times A\ni(x,y)\mapsto xy\in A$ on~$A$, called {\it multiplication} and satisfying
 the following for all $x,y,z\in A$ and ${\alpha}\in {\mathbb F}$:

 \subsec{(1)}~$\ xy=yx$;

 \subsec{(2)}~$\ (x+y)z=xz+yz$;

 \subsec{(3)}~$\ {\alpha}(xy)=({\alpha}x)y$;

 \subsec{(4)}~$\ ({x^2}y)x={x^2}(yx)$.

 An element $e$ of a Jordan algebra $A$ is a~{\it unit element}
 or a {\it unit} of~$A$, if
 $e\neq 0$ and $ea=a$ for all $a\in A$.

 \subsec{4.2.}~Recall that a $J\!B$-{\it algebra} is simultaneously
 a real Banach space~$A$ and a unital Jordan algebra with unit $\mathbb{1}$ such that

 \subsec{(1)}~$\ \|xy\|\le\|x\|\cdot\|y\| \quad (x, y\in A)$,

 \subsec{(2)}~$\ \|x^2\|=\|x\|^2 \quad (x\in A)$,

 \subsec{(3)}~$\ \|x^2\|\le\|x^2+y^2\| \quad(x, y\in A)$.

 The set $A_+:=\{x^2 : x\in A\}$, presenting a proper convex  cone, determines the structure of an ordered vector space on $A$ so that the unity $\mathbb{1}$ of the algebra $A$ serves as a strong  order unit, and the order interval $[-\mathbb{1},\,\mathbb{1}]:=\{x\in A : -\mathbb{1}\le x\le\mathbb{1}\}$ serves as the unit ball. Moreover, the inequalities $-\mathbb{1}\le x\le\mathbb{1}$ and
 $0\le x^2\le\mathbb{1}$ are equivalent.

 The intersection of all maximal associative subalgebras of $A$ is called the {\it center} of $A$ and denoted by $\mathcal{Z}(A)$.
 The element $a$ belongs to $\mathcal{Z}(A)$ if and only if $(ax)y=a(xy)$ for all $x, y\in A$. If $\mathcal{Z}(A)=\mathbb{R}\cdot\mathbb{1}$, then $A$ is said to be a $J\!B$-\textit{factor}. The center $Z(A)$ is an associative $J\!B$-algebra, and  such an algebra is isometrically isomorphic to the real Banach algebra $C(Q)$ of continuous functions on some compact space $Q$.

 \subsec{4.3.}~The idempotents of a $J\!B$-algebra are also called {\it projections}. The set of all projections $\mathbb{P}(A)$ forms a complete lattice with the order defined as $\pi\leq\rho\Longleftrightarrow\pi\circ\rho=\pi$. The sublattice of \textit{central projections\/} $\mathbb{P}_c(A)\!:=\mathbb{P}(A)\cap\mathcal{Z}(A)$ is a Boolean
 algebra. Assume that $\mathbb{B}$ is a subalgebra of the Boolean algebra $\mathbb{P}_c(A)$ or, equivalently,
 $\mathbb B\,(\mathbb R)$
 is a subalgebra of the center
 $\mathcal Z (A)$ of~$A$.
 Then we say that
 $A$ is
 a $\mathbb B$-$J\!B$-{\it algebra}
 if, for every partition of unity
 $(e_\xi)_{\xi \in \Xi}$
 in
 $\mathbb B$
 and every family
 $(x_\xi)_{\xi \in \Xi}$
 in
 $A$,
 there exists a unique $\mathbb B$-mixing
 $x:=\mix_{\xi \in \Xi}\,(e_\xi x_\xi)$,
 i.e., a unique element
 $x\in A$ such that
 $e_\xi x_\xi=e_\xi x$
 for all
 $\xi\in\Xi$.
 If
 $\mathbb B\,(\mathbb R)=\mathcal Z(A)$,
 then a $\mathbb B$-$J\!B$-algebra is also referred to as
 {\it centrally extended
 $J\!B$-algebra.}

 The unit ball of a $\mathbb{B}$-$J\!B$-algebra is closed under
 $\mathbb{B}$-mixings. Consequently, each $\mathbb{B}$-$J\!B$-algebra is a $\mathbb{B}$-cyclic Banach space.

 \subsec{4.4.}~\theorem{}The restricted descent of a $J\!B$-algebra in the model $\mathbb{V}^{(\mathbb{B})}$ is a $\mathbb{B}$-$J\!B$-algebra. Conversely, for every
 $\mathbb{B}$-$J\!B$-algebra $A$ there exists a unique $($up to isomorphism$)$ $J\!B$-algebra $\mathcal{A}$ within $\mathbb{V}^{\mathbb{B}}$ whose restricted descent is isometrically $\mathcal{B}$-isomorphic to $A$. Moreover,
 $[\![\mathcal{A}$ is a $J\!B$-factor $\!]\!]=\mathbb{1}$ if and only if\/ $\mathbb{B}\,(\mathbb{R})=\mathcal{Z}(A)$.
 \Endproc

 \beginproof~See \cite[Theorem 12.7.6]{IBA} and \cite[Theorem 3.1]{KusJB}.~\endproof

 \subsec{4.5.}~Now we give two applications of the above
 Boolean valued representation result to $\mathbb{B}$-$J\!B$-algebras. Theorems 4.7 and 4.11 below appear by transfer of the corresponding facts from the theory of $J\!B$-algebras.

 Let $A$ be a $\mathbb B$-$J\!B$-algebra and let $\Lambda :=\mathbb B\,(\mathbb R)$. An operator $\Phi\in A^\sharp $ is called a $\Lambda$-{\it valued state} if $\Phi\ge 0$ and $\Phi(\mathbb{1})=\mathbb{1}$. A state $\Phi$ is said to be {\it normal\/} if, for every increasing net $(x_{\alpha})$ in $A$ with the least upper bound $x:=\sup x_\alpha$, we have $\Phi(x)=\olim\Phi(x_\alpha)$. If $\mathcal{A}$ is the Boolean valued representation  of the algebra $A$, then the ascent $\varphi:=\Phi\!\uparrow$ is a bounded linear functional on $\mathcal{A}$ by
 \cite[Theorem 5.8.12]{BA_ST}. Moreover, $\varphi$ is positive and order continuous; i.e., $\varphi$ is a normal state on $\mathcal{A}$. The converse is also true: if $[\![\varphi$ is a  normal state on $\mathcal{A}]\!]=\mathbb{1}$, then the restriction of the operator $\varphi\!\downarrow$ to $A$ is a $\Lambda$-valued normal state. Now we will characterize $\mathbb{B}$-$J\!B$-algebras that are $\mathbb{B}$-dual spaces. Toward this end, it suffices to give Boolean valued interpretation for the following result.

 \subsec{4.6.}~\theorem{}A $J\!B$-algebra is a dual Banach space if
 and only if it is monotone complete and has a separating set of normal states.
 \Endproc

 \beginproof~See \cite[Theorem 2.3]{Shu}.~\endproof

 \subsec{4.7.}~\theorem{}Let $\mathbb{B}$ be a complete Boolean algebra and $\Lambda$ a Dedekind complete unital $AM$-space with $\mathbb{B}\simeq\mathbb{P}(\Lambda)$.
 A $\mathbb B$-$J\!B$-algebra $A$ is a $\mathbb{B}$-dual space if and only if $A$ is monotone complete and admits a separating set of $\Lambda$-valued normal states. If one of these equivalent conditions holds, then the part of $A^{\scriptscriptstyle\#}$
 consisting of order continuous operators serves as a $\mathbb{B}$-predual space of $A$. \Endproc

 \beginproof~See \cite[Theorem 12.8.5]{IBA} and \cite[Theorem 4.2]{KusJB}.~\endproof

 \subsec{4.8.}~An algebra $A$ satisfying one of the equivalent conditions 4.7 is called a $\mathbb B$-$JBW$-{\it algebra}. If, moreover, $\mathbb{B}$ coincides with the set
 of all central projections, then $A$ is said to be a $\mathbb{B}$-$JBW$-{\it factor}.
 It follows from Theorems 4.4 and 4.7 that $A$ is a
 $\mathbb B$-$JBW$-algebra $(\mathbb B$-$JBW$-factor) if and only if
 its Boolean valued representation $\mathcal A\in V^{(\mathbb B)}$ is a
 $JBW$-algebra ($JBW$-factor).

 A mapping $\Phi:A_+\to\Lambda\cup\{+\infty\}$ is a ($\Lambda$-valued) \textit{weight\/}
 if the conditions are satisfied (under the assumptions that $\lambda+(+\infty):=+\infty+\lambda:=+\infty$,
 $\lambda\cdot(+\infty)=:\lambda$ for all $\lambda\in\Lambda$, while
 $0\cdot(+\infty):=0$ and $+\infty+(+\infty):=+\infty$):

 \subsec{(1)}~$\Phi(x+y)=\phi(x)+\Phi(y)$ for all $x,y\in A_+$.

 \subsec{(2)}~$\Phi(\lambda x)=\lambda\Phi(x)$ for all $x\in A_+$ and $\lambda\in\Lambda_+$.

 A weight $\Phi$ is said to be a \textit{trace\/} if the additional condition is satisfied%

 \subsec{(3)}~$\Phi(x)=\Phi(U_sx)$ for all $x\in A_+$ and $s\in A$ with $s^2=\mathbb{1}$.

 Here, $U_a$ is the operator from $A$ to $A$ defined for a given
 $a\in A$ as $U_a:x\mapsto 2a(ax)-a^2$ $(x\in A)$. This operator is positive, i.e.,
 $U_a(A_+)\subset A_+$. If $a\in\mathcal{Z}(A)$, then $U_ax=a^2x$ $(x\in A)$.

 A weight (trace) $\Phi$ is called: \textit{normal\/} if $\Phi(x)=\sup_\alpha\Phi(x_\alpha)$
 for every increasing net $(x_\alpha)$ in $A_+$ with $x=\sup_\alpha x_\alpha$; \textit{semifinite\/}
 if there exists an increasing net $(a_\alpha)$ in $A_+$ with $\sup_\alpha a_\alpha=\mathbb{1}$
 and $\Phi(a_\alpha)\in\Lambda$ for all $\alpha$; \textit{bounded\/} if $\Phi(\mathbb{1})\in\Lambda$.
 Given two $\Lambda$-valued weights $\Phi$ and $\Psi$ on $A$, say that $\Phi$ is dominated by $\Psi$
 if there exists $\lambda\in\Lambda_+$ such that $\Phi(x)\leq\lambda\Psi(x)$ for all $x\in A_+$.

 \subsec{4.9.}~We need a few additional remarks about descents
 and ascents. Fix $+\infty\in\mathbb{V}^{(\mathbb{B})}$. If
 $\Lambda\!=\mathcal{R}{\Downarrow}$ and $\Lambda^\mathbb{u}\!=\mathcal{R}{\downarrow}$ then
 $$
 (\Lambda^\mathbb{u}\cup\{+\infty\}){\uparrow}=(\Lambda\cup\{+\infty\}){\uparrow}=
 \Lambda{\uparrow}\cup\{+\infty\}{\uparrow}=\mathcal{R}\cup\{+\infty\}.
 $$
 At the same time, $\Lambda^\star\!:=(\mathcal{R}\cup\{+\infty\}){\downarrow}=
 \mix(\mathcal{R}{\downarrow}\cup\{+\infty\})$
 consists of all elements of the form $\lambda_\pi\!:=\mix(\pi\lambda,\pi^\perp(+\infty))$
 with $\lambda\in\Lambda^{\mathbb{u}}$ and $\pi\in\mathbb{P}(\Lambda)$.
 Thus,  $\Lambda^\mathbb{u}\cup\{+\infty\}$ is a proper subset of
 $\Lambda^\star$, since $x_\pi\in\Lambda\cup\{+\infty\}$ if and only if
 $\pi=0$ or $\pi=I_\Lambda$.

 Assume now that $A=\mathcal{A}{\downarrow}$ with $\mathcal{A}$
 a $J\!B$-algebra within $\mathbb{V}^{(\mathbb{B})}$
 and $\mathbb{B}$ equal to $\mathbb{P}(A)$. Every bounded
 weight $\Phi:A\to\Lambda$ is evidently extensional:
 $b\!:=[\![x=y]\!]$ implies $bx=by$ which in turn yields
 $b\Phi(x)=\Phi(bx)=\Phi(by)=b\Phi(y)$ or, equivalently,
 $b\leq[\![\Phi(x)=\Phi(y)]\!]$. But an unbounded weight may fail to be
 extensional. Indeed, if $\Phi(x_0)=+\infty$ and $\Phi(x)\in\Lambda$
 for some $x_0,x\in A$ and $b\in\mathbb{P}(A)$ then
 $$
 \Phi(\mix(bx,b^\perp x_0))=
 \mix(b\Phi(x),b^\perp(+\infty))
 \notin\Lambda\cup\{+\infty\}.
 $$

 Given a semifinite weight $\Phi$ on $A$,
 we define its extensional modification $\widehat{\Phi}:A\to\Lambda^\star$ as follows:
 If $\Phi(x)\in\Lambda$ we put $\widehat{\Phi}(x)\!:=\Phi(x)$. If
 $\Phi(x)=+\infty$ then $x=\sup D$ with $D\!:=\{a\in A:\,0\leq a\leq x,\,\Phi(a)\in\Lambda\}$.
 Let $b$ stand for the greatest element of $\mathbb{P}(\Lambda)$ such that
 $\Phi(bD)$ is order bounded in $\Lambda^\mathbb{u}$ and put
 $\lambda\!:=\sup\Phi(bD)$. We define $\widehat{\Phi}(x)$ as $\lambda_b=
 \mix(b\lambda,b^\perp(+\infty))$; i.e., $b\widehat{\Phi}(x)=\lambda$ and
 $b^\perp\widehat{\Phi}(x)=b^\perp(+\infty)$. It is not difficult to
 check that $\widehat{\Phi}$ is an extensional mapping. Thus, for
 $\varphi\!:=\widehat{\Phi}{\uparrow}$ we have
 $[\![\varphi:\mathcal{A}\to\mathcal{R}\cup\{+\infty\}]\!]=\mathbb{1}$
 and, according to \cite[1.6.6]{BA_ST}, $\widehat{\Phi}=\varphi{\downarrow}\ne\Phi$.
 But if we define $\varphi{\Downarrow}$ as $\varphi{\Downarrow}(x)=\varphi{\downarrow}(x)$
 whenever $\varphi{\downarrow}(x)\in\Lambda$ and
 $\varphi{\Downarrow}(x)=+\infty$ otherwise, then
 $\Phi=(\widehat{\Phi}{\uparrow}){\Downarrow}$.

 \subsec{4.10.}\theorem{}Let $A$ be a $J\!BW$-algebra and let $\tau$ be a normal semifinite real-valued trace on $A$. For each real-valued weight $\varphi$ on $A$ dominated by $\tau$ there exists a unique positive element $h\in A$ such that $\varphi(a)=\tau(U_{h^{1/2}}a)$ for all $a\in A_+$. Moreover, $\varphi$ is bounded if and only if $\tau(h)$ is finite and $\varphi$ is a trace if and only if $h$ is a central element of $A$.
 \Endproc

 \beginproof~This fact was proved in \cite{King}.~\endproof

 \subsec{4.11.}~\theorem{}Let $A$ be a $\mathbb{B}$-$J\!BW$-algebra and $\mathrm{T}$ is a normal semifinite $\Lambda$-valued trace on $A$. For each weight $\Phi$ on $A$ dominated by $\mathrm{T}$ there exists a unique positive $h\in A$ such that $\Phi(x)= \mathrm{T}(U_{h^{1/2}}x)$ for all $x\in A_+$. Moreover, $\Phi$ is bounded if and only if $\mathrm{T}(h)\in\Lambda$ and $\Phi$ is a trace if and only if $h$ is a central element of $A$.
 \Endproc

 \beginproof~We present a sketch of the proof. Taking into consideration the remarks in 4.9, we define $\varphi\!=\widehat{\Phi}{\uparrow}$ and $\psi\!=\widehat{\Psi}{\uparrow}$. Then within $\mathbb{V}^{(\mathbb{B})}$ the following hold: $\tau$ is a semifinite normal real-valued trace on $\mathcal{A}$ and $\varphi$ is real-valued weight on $\mathcal{A}$ dominated by $\tau$. By transfer principle we may apply Theorem 4.10 and find $h\in\mathcal{A}$ such that $\varphi(x)=\tau(U_{h^{1/2}}x)$ for all $x\in\mathcal{A}_+$. Actually, $h\in A$ and $\varphi{\Downarrow}(x)=\tau{\Downarrow}(U_{h^{1/2}}x)$ for all $x\in A_+$. It remains to note that $\Phi=\varphi{\Downarrow}$ and $\mathrm{T}=\tau{\Downarrow}$. The details of the proof are left to the reader.~\endproof

 \subsec{4.12.}~~$J\!B$-algebras are nonassociative
 real analogs of $C^*$-algebras and von Neumann operator algebras.
 The theory of these algebras stems from  Jordan, von~Neumann, and Wigner \cite{JNW}  and exists as a~branch of functional analysis since the mid 1960s. The stages of its development are reflected in Alfsen, Shultz, and
 St\o rmer \cite{ASS}. The theory of $J\!B$-algebras undergoes intensive study, and the scope of its applications widens.
 Among the main  areas of research are the structure and classification of $J\!B$-algebras, nonassociative integration and
 quantum probability theory, the geometry of states of $J\!B$-algebras, etc.; see Ajupov \cite{Aju1, Aju2}; Hanshe-Olsen and St\"ormer \cite{H-OS} as well as the references therein.

 {\bf(2)}~The Boolean valued approach to $JB$-algebras was charted by
 Kusraev in~the article \cite{KusJB} which contains Theorems 4.4 and 4.7 (also see \cite{KusJB1}).
 These theorems are instances of the Boolean valued interpretation of the results
 by Shulz \cite{Shu} and  by Ajupov and Abdullaev \cite{AjA}. In \cite{KusJB} Kusraev introduced the class of $\mathbb{B}$-$JBW$-algebras which is broader than the class of $JBW$-algebras. The principal distinction is that
 a~$\mathbb{B}$-$JBW$-algebra has faithful representation as the algebra of selfadjoint operators in some $AW^ \ast$-module rather than on a~Hilbert space as in the case of $JBW$-algebras (cp.~Kusraev and Kutateladze~\cite{IBA}). The class of $AJW$-algebras was firstly mentioned by Topping in~\cite{Top}. Theorem 4.11 was not published before.

 \section*{5.~Transfer in Harmonic Analysis}

 In what follows, $G$ is a locally compact abelian group and $\tau$ is the topology of $G$, while $\tau(0)$ is a~neighborhood base of $0$ in $G$ and $G'$ stands for the dual group of~$G$. Note that $G$ is also the dual group of $G'$ and we write $\langle x,\gamma\rangle\!:=\gamma(x)$ $(x\in G,\,\gamma\in G')$.

 \subsec{5.1.}~By restricted transfer,  $G^{\scriptscriptstyle\wedge}$ is a group within $\mathbb{V}^{(\mathbb{B})}$. At the same time $\tau(0)^{\scriptscriptstyle\wedge}$ may fail to be a topology on $G^{\scriptscriptstyle\wedge}$. But $G^{\scriptscriptstyle\wedge}$ becomes a~topological group on taking $\tau(0)^{\scriptscriptstyle\wedge}$ as a neighborhood base of $0\!:=0^{\scriptscriptstyle\wedge}$. This topological group is again denoted by $G^{\scriptscriptstyle\wedge}$ itself. Clearly, $G^{\scriptscriptstyle\wedge}$ may not be locally
 compact. Let $U$ be a neighborhood of $0$ such that $U$ is compact. Then $U$ is totally bounded.
 It follows  by restricted transfer  that $U^{\scriptscriptstyle\wedge}$ is  totally bounded as well, since total boundedness can be expressed by a~restricted formula. Therefore the completion of $G^{\scriptscriptstyle\wedge}$ is  locally compact. The completion of $G^{\scriptscriptstyle\wedge}$ is denoted by $\mathcal{G}$, and by the above observation $\mathcal{G}$ is a locally compact abelian group within $\mathbb{V}^{(\mathbb{B})}$.

 \subsec{5.2.}~Let $Y$ be a Dedekind complete vector lattice and let $Y_\mathbb{C}$ be the complexification
 of~$Y$. A vector-function $\varphi:G\to Y$ is said to be {\it uniformly order continuous\/} on a~set $K$ if
 $$
 \inf\limits_{U\in\tau(0)}\sup\{|\varphi(g_1)-\varphi (g_2)|: \ g_1,g_2\in K,\ g_1-g_2\in U\}=0.
 $$
 This amounts to saying that $\varphi$ is order bounded on $K$ and, if $e\in Y$ is an arbitrary upper bound of $\varphi(K)$, then for each $0<\varepsilon\in\mathbb{R}$ there exists a partition of unity $(\pi_\alpha)_{\alpha\in\tau(0)}$ in $\mathbb{P}(Y)$ such that $\pi_\alpha|\varphi(g_1)-\varphi (g_2)|\leq\varepsilon e$ for all $\alpha\in\tau(0)$ and $g_1,g_2\in K$, $g_1-g_2\in\alpha$. If, in this definition we put $g_2=0$, then we arrive at the definition of mapping  \textit{order continuous at zero}.

 We now introduce the class of dominated $Y$-valued mappings.
 A~mapping $\psi: G\to Y_\mathbb{C}$ is called {\it positive definite\/} if
 $$
 \sum\limits_{j,k=1}^n \psi(g_j-g_k)\ c_j \overline{c}_k\geq0
 $$
 for all finite collections $g_1,\ldots,g_n \in G$ and
 $c_1,\ldots,c_n\in\mathbb{C}$ $(n\in\mathbb{N})$.

 For $n=1$, the definition  implies readily that $\psi(0)\in Y_+$. For
 $n=2$, we have  $|\psi(g)|\le\psi(0)$ $(g\in G)$.
 If we introduce the structure of an $f$-algebra with unit $\psi(0)$
 in the order ideal of $Y$ generated by $\psi(0)$ then, for $n=3$,
 from the above definition we can deduce one more inequality
 $$
 |\psi(g_1)-\psi(g_2)|^2\le 2 \psi(0)(\psi(0)-\Re\psi(g_1-g_2)) \quad
 (g_1,g_2\in G).
 $$
 It follows that every positive definite mapping $\psi : G\to Y_\mathbb{C}$
 $o$-continuous at zero is order-bounded (by the element $\psi(0)$)
 and uniformly $o$-continuous.
 A mapping $\varphi : G\to Y$ is called {\it dominated\/} if there exists
 a~positive definite mapping $\psi : G\to Y_\mathbb C$ such that
 $$
 \Bigg|\sum\limits_{j,k=1}^n \varphi(g_j-g_k) c_j
 \overline{c}_k\Bigg|\le \sum\limits_{j,k=1}^{n}\psi(g_j-g_k) c_j \overline
 c_k
 $$
 for all $g_1,\ldots,g_n \in G$, $c_1,\ldots, c_n \in\mathbb{C}$ $(n\in\mathbb{N})$. In this case we also say that $\psi$ is a {\it dominant} of $\varphi$. It can be easily shown that if $\varphi:G\to Y_\mathbb{C}$ has  dominant order continuous at zero then $\varphi$ is order bounded and uniformly order continuous.

 We denote by $\mathfrak{D}(G,Y_\mathbb{C})$ the vector space of all dominated mappings from $G$ into $Y_\mathbb{C}$ whose dominants are order continuous at zero. We also consider the set $\mathfrak{D}(G,Y_\mathbb{C})_+$ of all positive definite mappings from $G$ into $Y_\mathbb{C}$. This set is a proper cone in $\mathfrak{D}(G,Y_\mathbb{C})$ and defines the order  compatible with the structure of a~vector space on $\mathfrak{D}(G,Y_\mathbb{C})$. Actually, $\mathfrak{D}(G,Y_\mathbb{C})$ is a Dedekind complete complex vector lattice; cp.~5.13 below. Also, define
 $\mathfrak{D}(\mathcal{G},\mathcal{C})\in\mathbb{V}^{(\mathbb{B})}$ to be the set of functions $\varphi:\mathcal{G}\to\mathcal{C}$ with the property that $[\![\varphi$ has  dominant  continuous at zero$]\!]=\mathbb{1}$.

  \subsec{5.3.}\proclaim{}Let $Y=\mathcal{R}{\downarrow}$. For every
 $\varphi\in\mathfrak{D}(G,Y_\mathbb{C})$ there exists a unique $\tilde{\varphi}\in\mathbb{V}^{(\mathbb{B})}$
 such that $[\![\tilde{\varphi}\in\mathfrak{D}(\mathcal{G},\mathcal{C})]\!]=\mathbb{1}$
 and
 $[\![\tilde{\varphi}(x^{\scriptscriptstyle\wedge})=\varphi(x)]\!]=\mathbb{1}$
 for all $x\in G$. The mapping $\varphi\mapsto\tilde{\varphi}$ is an
 linear and order isomorphism from $\mathfrak{D}(G,Y)$ onto
 $\mathfrak{D}(\mathcal{G},\mathcal{C}){\downarrow}$.
 \Endproc

 \subsec{5.4.}~Define $C_0(G)$ as the space of all
 continuous complex functions $f$ on $G$ vanishing at infinity. The latter
 means that for every $0<\varepsilon\in\mathbb{R}$ there exists a
 compact set $K\subset G$ such that $|f(x)|<\varepsilon$ for all $x\in G\setminus
 K$. Denote by $C_{c}(G)$ the space of all
 continuous complex functions on $G$ having compact support. Evidently, $C_c(G)$
 is dense in $C_0(G)$ with respect to the norm $\|\cdot\|_\infty$.

 \subsec{5.5.}~Introduce the class of dominated operators. Let $X$ be a complex normed space and let $Y$ be a complex Banach lattice. A linear operator $T:X\to Y$ is said to be   \textit{dominated\/} or
 {\it having abstract norm\/}if $T$ sends the unit ball of $X$ into an order bounded subset of $Y$. This amounts to saying that there exists $c\in Y_+$ such that $|Tx|\leq c\|x\|_\infty$ for all $x\in C_0(Q)$. The~set of all dominated operators from $X$ to $Y$
 is denoted by $L_m(X, F)$.  If $Y$ is Dedekind complete then
 $$
 \[T\]\!:=\{|Tx|:\ x\in X,\ \|x\|\leq 1\}
 $$
 exists and is called the~{\it abstract norm\/} or {\it dominant\/} of $T$. Moreover, if $X$ is a vector lattice and $Y$ is Dedekind complete then $L_m(X,Y)$ is a vector sublattice of $L^{\sim}(X,Y)$.

 Given a positive   $T\in L_m(C_0(G'),Y)$, we can define the mapping $\varphi:G\to Y$ by putting $\varphi(x)=T(\langle x,\cdot\rangle)$ for all $(x\in G)$, since  $\gamma\mapsto\langle x,\gamma\rangle$ lies in $C_0(G')$ for every $x\in G$. It is not difficult to ensure that   $\varphi$ is order continuous at zero and positive definite. The converse is also true; see 5.8.

 \subsec{5.6.}~Consider a metric space $(M,r)$.
 The definition of  metric space can be written as a bounded formula, say $\varphi(M,r,\mathbb{R})$, so that $[\![\varphi(M^{\scriptscriptstyle\wedge},
 r^{\scriptscriptstyle\wedge},
 \mathbb{R}^{\scriptscriptstyle\wedge})]\!]=\mathbb{1}$ by restricted transfer.
 In other words, $(M^{\scriptscriptstyle\wedge},
 r^{\scriptscriptstyle\wedge})$ is a metric space within $\mathbb{V}^{(\mathbb{B})}$. Moreover we consider the internal function $r^{\scriptscriptstyle\wedge}:\, M^{\scriptscriptstyle\wedge}\to
 \mathbb{R}^{\scriptscriptstyle\wedge}\subset\mathcal{R}$ as an $\mathcal{R}$-valued metric on $M^{\scriptscriptstyle\wedge}$.
 Denote by $(\mathcal{M},\rho)$ the completion of   $(M^{\scriptscriptstyle\wedge}, r^{\scriptscriptstyle\wedge})$; i.e., $[\![(\mathcal{M},\rho)$ is a complete metric space $]\!]=\mathbb{1}$,
 $[\![M^{\scriptscriptstyle\wedge}$ is a dense subset of $\mathcal{M}]\!]=\mathbb{1}$, and $[\![r(x^{\scriptscriptstyle\wedge})=
 \rho(x^{\scriptscriptstyle\wedge})]\!]=\mathbb{1}$ for all $(x\in M)$

 Now, if $(X,\|\cdot\|)$ is a real (or complex) normed space then  $[\![X^{\scriptscriptstyle\wedge}$ is a vector space over the field $\mathbb{R}^{\scriptscriptstyle\wedge}$ (or $\mathbb{C}^{\scriptscriptstyle\wedge}$) and $\|\cdot\|^{\scriptscriptstyle\wedge}$ is a norm on
 $X^{\scriptscriptstyle\wedge}$ with values in
 $\mathbb{R}^{\scriptscriptstyle\wedge}\subset\mathcal{R}
 ]\!]=\mathbb{1}$. So, we will consider $X^{\scriptscriptstyle\wedge}$ as
 an~$\mathbb{R}^{\scriptscriptstyle\wedge}$-vector space with  $\mathcal{R}$-valued norm within $\mathbb{V}^{(\mathbb{B})}$.
 Let $\mathcal X\in\mathbb{V}^{(\mathbb{B})}$ stand for the (metric) completion of $X^{\scriptscriptstyle\wedge}$ within $\mathbb{V}^{(\mathbb{B})}$. It is not difficult to see that
 $[\![\mathcal{X}$ is a real (complex) Banach space including $X^{\scriptscriptstyle\wedge}$ as an
 $\mathbb{R}^{\scriptscriptstyle\wedge}(\mathbb{C}^{\scriptscriptstyle\wedge})$-linear
 subspace$]\!]=\mathbb{1}$, since the metric $(x,y)\mapsto\|x-y\|$ on
 $X^{\scriptscriptstyle\wedge}$ is translation invariant.
 Clearly, if $X$ is a real (complex) Banach lattice then
 $[\![\mathcal{X}$ is a real (complex) Banach lattice including $X^{\scriptscriptstyle\wedge}$
 as an $\mathbb{R}^{\scriptscriptstyle\wedge}(\mathbb{C}^{\scriptscriptstyle\wedge})$-linear
 sublattice$]\!]=\mathbb{1}$.

 \subsec{5.7.}~\Theorem{}Let $Y=\mathcal{C}{\downarrow}$ and $\mathcal{X}'$ be the topological dual of $\mathcal{X}$ within $\mathbb{V}^{(\mathbb{B})}$. For every $T\in L_m(X,Y)$ there exists a unique $\tau\in\mathcal{X}'{\downarrow}$ such that $[\![\tau(x^{\scriptscriptstyle\wedge})=T(x)]\!]=\mathbb{1}$ for all $x\in X$. The mapping $T\mapsto\phi(T)\!:=\tau$ defines an isomorphism between the $\mathcal{C}{\downarrow}$-modules $L_m(X,Y)$ and $\mathcal{X}'{\downarrow}$. Moreover, $\[T\]=\[\phi(T)\]$ for all $T\in L_m(X,Y)$. If $X$ is a normed lattice then $[\![\mathcal{X}'$ is a Banach lattice $]\!]=\mathbb{1}$, while $\mathcal{X}'{\downarrow}$ is a vector lattice and $\phi$ is a lattice isomorphism.
 \Endproc

 \beginproof~Suffice it to consider the real case. Apply  \cite[Theorem 8.3.2]{DOP} to the lattice normed space $X\!:=(X,\[ \cdot\])$ with $\[ x\]=\|x\|\mathbb{1}$. By \cite[Theorem 8.3.4\,(1) and Proposition 8.3.4\,(2)]{DOP} the spaces
 $\mathcal{X}'{\downarrow}:=\mathcal{L}^{(\mathbb{B})}(\mathcal{X}, \mathcal{R}){\downarrow}$ and $L_m(X,Y)$ are linear isometric. We are left with referring to \cite[Proposition 5.5.1\,(1)]{DOP}.~\endproof

 \subsec{5.8.}~\Theorem{}A mapping $\varphi:G\to Y_\mathbb{C}$ is order continuous at zero and positive definite if and only if there exists a unique positive operator $T\in L_m(C_0(G'),Y_\mathbb{C})$ such that $\varphi(x)=T(\langle x,\cdot\rangle)$ for all $(x\in G)$.
 \Endproc

 \beginproof~By transfer, 5.3, and Theorem 5.7, we can replace $\varphi$ and $T$ by their Boolean valued representations $\tilde{\varphi}$ and $\tau$. The norm completion of $C_0(G')^{\scriptscriptstyle\wedge}$ within $\mathbb{V}^{(\mathbb{B})}$ coincides with $C_0(\mathcal{G}')$. (This can be proved by the reasoning similar to that in Takeuti\author{Takeuti G.} \cite[Proposition 3.2]{Tak2}.) Application of the classical Bochner Theorem (see Loomis\author{Loomis L.H} \cite[Section 36A]{Loo}) to $\tilde{\varphi}$ and $\tau$ yields the desired result.~\endproof

 \subsec{5.9.}~We now specify the vector integral of use  in this
  subsection; see details in~\cite[5.14.B]{BA_ST}.
  Let $\mathcal{A}$ be a $\sigma$-algebra of subsets of $Q$, i.e. $\mathcal{A}\subset\mathcal{P}(Q)$. We identify  this algebra with the isomorphic algebra of the characteristic functions $\{1_A\!:=\chi_A:\,A\in\mathcal{A}\}$ so that $S(\mathcal A)$ is the space of all $\mathcal{A}$-simple functions on $Q$; i.e., $f\in S(\mathcal{A})$ means that $f=\sum_{k=1}^n\alpha_k\chi_{A_k}$ for some $\alpha_1,\dots,\alpha _n\in\mathbb{R}$ and 
 disjoint $A_1,\dots,A_n\in\mathcal{A}$. Let a measure $\mu$ be defined on $\mathcal{A}$ and take values in a~Dedekind complete vector lattice~$Y$. We suppose that $\mu$ is order bounded. If $f\in S(\mathcal{A})$ then we put
 $$
 I_\mu \!:=\int f\, d\mu =\sum\limits_{k=1}^n\alpha_k\mu(A_k).
 $$
 It can be easily seen that the integral $I_\mu $ can be extended to the spaces of $\mu$-summable functions
 $\mathcal{L}^1(\mu)$ for which the more informative notations $\mathcal{L}^1(Q,\mu)$ and $\mathcal{L}^1(Q,\mathcal A,\mu)$ are also used. On identifying equivalent functions, we obtain the Dedekind $\sigma$-complete vector lattice $L^1(\mu)\!:=L^1(Q,\mu)\!:=L^1(Q,\mathcal A,\mu)$.

 \subsec{5.10.}~Assume now that $Q$ is a topological space. Denote by $\mathcal{F}(Q)$, $\mathcal{K}(Q)$, and $\mathcal{B}(Q)$  the collections of all closed, compact, and Borel subsets of $Q$. A measure $\mu:\mathcal{B}(Q)\to Y$ is said to be {\it quasi-Radon} (\textit{quasi-regular}) if $\mu$ is order bounded and
 $$
 \gathered
 |\mu|(U)=\sup\{|\mu|(K):\ K\in\mathcal{K}(Q),\,K\subset U\}
 \\
 (|\mu|(U)=\sup\{|\mu|(K):\ K\in\mathcal{F}(Q),\,K\subset U\}).
 \endgathered
 $$
 for every open (respectively, closed) set $U\subset Q$. If the above  equalities are fulfilled for all Borel $U\subset Q$ then we speak about {\it Radon} and \textit{regular\/} measures. Say that $\mu=\mu_1+i\mu_2:\mathcal{B}(Q)\to Y_\mathbb{C}$ have one of the above properties whenever this property
 is enjoyed by both $\mu_1$ and $\mu_2$. We denote by $\qca(Q,Y)$ the vector lattice of all $\sigma$-additive  $Y_\mathbb{C}$-valued quasi-Radon measures on $\mathcal{B}(Q)$. If $Q$ is locally compact or even completely regular then $\qca(Q,Y)$ is a vector lattice; see \cite[Theorem 6.2.2]{DOP}.
 The variation of a $Y_{\mathbb{C}}$-valued (in particular,
 $\mathbb{C}$-valued) measure $\nu $ is denoted in the standard fashion: $|\nu|$.

 \subsec{5.11.}~\Theorem{}Let $Y$ be a real Dedekind complete vector lattice and let $Q$ be a locally compact topological space. Then for each  $T:L_m(C_0(Q),Y_\mathbb{C})$ there exists a~unique measure $\mu\!:=\mu_T\in\qca(Q,Y_\mathbb{C})$ such that
 $$
 T(f)=\int\limits_{Q}f\,d\mu
 \quad (f \in C_0(Q)).
 $$
 The mapping $T\mapsto\mu_T$ is a lattice isomorphism from $L_m(C_0(Q),Y_\mathbb{C})$ onto $\qca(Q,Y_\mathbb{C})$.
 \Endproc

 \beginproof~See \cite[Theorem 2.5]{KM5}.~\endproof

 \subsec{5.12.}~\Theorem{}Let $G$ be a locally compact abelian
 group, let $G'$ be the dual group of~$G$, and let $Y$ be a Dedekind complete real vector lattice.
 For a mapping $\varphi:G\to Y_\mathbb{C}$
 the following  are equivalent:

 \subsec{(1)}~$\varphi$ has  dominant  order continuous at zero.

 \subsec{(2)}~There exists a unique measure $\mu\in\qca(G',Y_\mathbb{C})$ such
 that
 $$
 \varphi(g)=\int\limits_{G'}\chi(g)\,d\mu(\chi)\quad(g\in G).
 $$
 \Endproc

 \beginproof~This is immediate from Theorems 5.8 and 5.11.~\endproof

 \subsec{5.13.}~\Corollary{}The Fourier transform establishes an order and linear isomorphism between the space
 of measures $\qca(G',Y)$  and the space of dominated mappings $\mathfrak{D}(G,Y_\mathbb{C})$. In particular, $\mathfrak{D}(G,Y_\mathbb{C})$ is a~Dedekind complete complex vector lattice.
 \Endproc

 \subsec{5.14.}~{\bf(1)}~In~\cite{Tak2} Takeuti introduced
 the Fourier transform for the mappings defined on a~locally compact
 abelian group and having as values  pairwise commuting normal
 operators in a~Hilbert space. By applying the transfer principle,
 he developed a~general technique for translating
 the classical results to operator-valued functions. In particular, he
  established a~version of the Bochner Theorem describing
 the set of all inverse Fourier transforms of positive
 operator-valued Radon measures. Given a~complete Boolean algebra
 $\mathbb{B}$ of projections in a~Hilbert space $H$, denote by
 $(\mathbb{B})$ the space of all selfadjoint operators on $H$ whose
 spectral resolutions are in $\mathbb{B}$; i.e., $A\in(\mathbb{B})$ if and only if $A=\int_\mathbb{R}\lambda\,dE_\lambda$ and $E_\lambda\in\mathbb{B}$ for all $\lambda\in\mathbb{R}$. If $Y\!:=(\mathbb{B})$ then Theorem 5.8 is essentially Takeuti's result \cite[Theorem~1.3]{Tak2}.

 {\bf(2)}~Kusraev and Malyugin in~\cite{KM5} abstracted  Takeuti's results in the
 following directions: First, they considered more general arrival spaces, namely,
 norm complete lattice normed spaces. So the important particular cases of Banach spaces and Dedekind complete
 vector lattices were covered. Second, the class of dominated mappings was identified with the set of all inverse Fourier
 transforms of order bounded quasi-Radon vector measures. Third, the construction of a suitable Boolean valued universe was eliminated from all definitions and statements of results.

 In particular, Theorem 5.12 and Corollary 5.13 correspond
 to \cite[Theorem 4.3]{KM5} and \cite[Theorem 4.4]{KM5}; while their
 lattice normed versions, to \cite[Theorem 4.1]{KM5} and \cite[Theorem 4.5]{KM5}, respectively.

 {\bf(3)}~Theorem 5.7 is due to Gordon \cite[Theorem 2]{Gor2}. Proposition~3.3 in Takeuti \cite{Tak2} is essentially the same result for the particular departure  and  arrival spaces; i.e., $X=L^1(G)$ and $Y=(\mathbb{B})$.
 Theorem 5.11 is taken from Kusraev and Malyugin \cite{KM5}. In the case of $Q$ compact, it was proved by
 Wright in~\cite[Theorem 4.1]{Wr6}. In this result $\mu$ cannot be chosen regular rather than quasiregular.
 The quasiregular measures were introduced by Wright in~\cite{Wr1}.

\bibliographystyle{plain}

\medskip

\noindent
{\it Anatoly G.~Kusraev}\\
{\leftskip\parindent\small
\noindent
Southern Mathematical Institute\\
 22 Markus Street\\
Vladikavkaz, 362027, RUSSIA\\
E-mail: kusraev@smath.ru
\par}

\medskip
\noindent
{\it Sem\"en S.~Kutateladze}\\
{\leftskip\parindent\small
\noindent
Sobolev Institute of Mathematics\\
4 Koptyug Avenue\\
Novosibirsk, 630090, RUSSIA\\
E-mail: sskut@math.nsc.ru
\par}
\medskip

\end{document}